\documentstyle[12pt]{article}
\newfont{\frak}{eufm10 scaled\magstep1}
\newfont{\extra}{msbm10 scaled\magstep1}

\newcommand{\sect}[1]{\setcounter{equation}{0}\section{#1}}
\newcommand{\subsect}[1]{\subsection{#1}}

\newcommand{\R}{\mbox{$I \!\! R \, $}}

\newcommand{\ba}{\begin{array}}
\newcommand{\ea}{\end{array}}

\newcommand{\be}{\begin{equation}}
\newcommand{\ee}{\end{equation}}
\newcommand{\eea}{\end{eqnarray}}
\newcommand{\bea}{\begin{eqnarray}}

\newcommand{\G}{{G}}
\newcommand{\g}{{g}}          
\begin{document}
\begin{center} 
{\LARGE{\bf{Contractions on the Classical Double}}}
\end{center}  
\bigskip
 
\begin{center} 
Angel Ballesteros$^1$ and Mariano A. del Olmo$^2$
\end{center}

\begin{center}
{\it{$^1$ Departamento de F\'{\i}sica,\\
Universidad de Burgos. E-09001, Burgos, Spain \\ \medskip
$^2$ Departamento de F\'{\i}sica Te\'orica,\\
Universidad de Valladolid. E-47011, Valladolid, Spain.}}
\end{center}

\vskip 0.5cm

\begin{abstract}

Lie algebra contractions on the classical Drinfel'd Double of a given Lie
bialgebra are introduced and compared to the usual Lie bialgebra
contraction theory. The connection between both approaches turns out to be
intimately linked to duality problems. The non-relativistic (Galilean)
limit of a (1+1) Poincar\'e Double is used to illustrate the contraction
process. Finally, it is shown that, in a certain sense, the classical
limit in a quantum algebra can be thought as a certain contraction on the
corresponding Double. 

\end{abstract}

\vskip 1cm


\sect{Introduction}

The well known result of Drinfel'd \cite{Dr} that establishes an one to
one correspondence between Poisson-Lie groups and Lie bialgebras is
the key stone to understand the essential role that Lie bialgebras
play in the quantization of Poisson-Lie structures (quantum groups).
The concept of classical double is just a reformulation of that of  Lie 
bialgebra in terms of a (double) dimensional Lie algebra. This
``duplication" process can be iterated by taking into account that
the double of a Lie bialgebra 
can be in turn equipped with a Lie bialgebra structure 
by means of a canonical $r$-matrix. The quantization of such a double
Lie bialgebra is the so-called quantum double, a structure that
has been essential in the explicit obtention of many quantum $R$-matrices.
(See \cite{CP} for a detailed exposition and references therein).  

On the  other hand, the initial developements in quantum groups were  
mainly devoted to the construction of deformations of semisimple
Lie structures \cite{FRT}, and it was soon discovered that
contraction methods allowed to construct quantum deformations
of  non-semisimple algebras in a very efficient way \cite{BCGS}. Since
then, contraction techniques have been extensively used to generate
quantum deformations of Poincar\'e and Galilei groups and algebras
(see, for instance, \cite{CGe2,Lukap,BHOS4d}). 

More recently, a systematic approach to the contractions of
quantum groups and algebras has been introduced \cite{BGHOS,BT,BO}.
That scheme is based in a contraction theory for the
underlying Lie bialgebras, that turn out to be again the
objects characterizing the behaviour of the full quantum structures.

The aim of this paper is to explore how contractions can
be implemented onto the corresponding classical doubles and to  
show that the resultant approach generalizes the previuos ones and allows
the study of new processes as the classical limit in a contraction
context. In Section 2 we shall briefly
review the notion of the classical double and in Section 3 the notions
of Lie algebra and Lie bialgebra  contraction will be recalled.
General contractions on the double Lie algebra 
$D$ of an arbitrary Lie bialgebra $(\g,\eta)$ 
will be characterized in Section 4.                    
In this respect it will be essential to consider how 
the internal pairing among
generators of the double behaves under contraction. The solution
to this problem will lead
to us to recover the Lie bialgebra contraction theory. 
Throughout the  paper, the approach will be illustrated 
by studying how the classical double of a Lie bialgebra
of the two-dimensional Euclidean algebra can give rise to a 
(1+1) Galilei double by means of a suitable contraction
that implements the non-relativistic limit.

\sect{The classical double of a Lie bialgebra}

The notion of Lie bialgebra $(\g,\eta)$ arises when  a Lie algebra
$\g$ is furnished with a cocommutator $\eta:\g\to \g\otimes \g$ such that
$\eta$ is a 1--cocycle and the dual map $\eta^{*}:\g^{*}\otimes \g^{*} \to
\g^{*}$ is a Lie bracket on $\g^*$ (the dual vector space of $\g$). A Lie
bialgebra
$(\g,\eta)$ is called  a coboundary bialgebra if there exists an element
$\rho\in
\g\otimes \g$ (the classical $r$--matrix), such that
\be
\eta(X)=[1\otimes X + X \otimes 1,\,  \rho], \qquad  \forall X\in
\g. \label{ag}
\ee

On the other hand, the map (\ref{ag}) defined by using an arbitrary
$\rho \in \g\otimes \g$ defines a Lie bialgebra if and only if the
symmetric  part of $\rho$ is $\g$--invariant and the antisymmetric part of
$\rho$ is a solution of the generalized Classical Yang--Baxter Equation
(CYBE) \cite{Dr}.

Let ($\cal G,\eta$) be a Lie bialgebra, and $\{X_i\}$ a basis of $\g$.
With such a basis, the Lie bialgebra can be characterized
 by the pair of structure tensors
($c^k_{ij},f^{lm}_n$) that define the commutator and the cocommutator
in the form:
\be
[X_i,X_j]=c^k_{ij}X_k, \quad \eta(X_n)= f^{lm}_n X_l\otimes X_m.  
\label{aga}\ee
In this language, the cocycle condition 
becomes the following compatibility condition between the
tensors $c$ and $f$
\be
f^{ab}_k c^k_{ij} = f^{ak}_i c^b_{kj}+f^{kb}_i c^a_{kj}
+f^{ak}_j c^b_{ik} +f^{kb}_j c^a_{ik}. \label{agb}
\ee  
Let us now consider a basis $\{x^i\}$ of $\g ^*$ such that 
$\langle x^i,X_j\rangle= \delta^i_j$;
 then ($\cal G ^*,\xi$) is also a Lie
bialgebra with structure tensors ($f, c$), i.e.,
\be
\{x^i,x^j\}=f^{ij}_k x^k, \quad \xi(x^n)=c^{n}_{lm} x^l\otimes x^m.  
\label{agc}
\ee

This intrinsic self-dual character between the two Lie algebras
included within a Lie bialgebra lead to the consideration of
the pair ($\cal G ,\cal G ^*$) and its associated vector space 
$\cal G \oplus \cal G ^*$, that can be endowed with a Lie algebra
structure by means of the commutators
\be
[X_i,X_j]= c^k_{ij}X_k, \quad  
[x^i,x^j]\equiv \{x^i,x^j\}= f^{ij}_k x^k, \quad
[x^i,X_j]= c^i_{jk}x^k- f^{ik}_j X_k.\label{agd}
\ee
This Lie algebra, $\cal D ({\cal G})$, is called the
Double Lie algebra of $(G,\eta)$. Obviously,
$\cal G$ and $\cal G^*$ are subalgebras of $\cal D$, and the 
compatibility condition (\ref{agb}) is just the Jacobi identity
for (\ref{agd}).

Moreover, if $\cal G$ is a finite dimensional Lie algebra,
then $\cal D$  is a
(coboundary) Lie bialgebra with classical $r$-matrix
\be
r=X_i\otimes x^i,
\label{rmat}
\ee
fulfilling the Classical Yang-Baxter Equation (CYBE).
The cocommutator $\delta(r)$  derived from (\ref{rmat}) is
\be
\delta(X_i)=f_i^{jk}\,X_j\wedge X_k,\qquad 
\delta(x^i)=c^i_{jk}\,x^j\wedge x^k.
\label{codob}
\ee
To summarize, we can say that $\cal D$ can be seen as 
both a Lie algebra and a (triangular) coboundary Lie bialgebra.
In fact this ``double Lie bialgebra" has as sub-Lie-bialgebras the 
original one $(\cal G,\eta)$ and its dual. This duplication process
can be obviously iterated.
It is also worthy to recall the following inner
product on $\cal G \oplus \cal G ^*$ \cite{CP}:
\be
 \langle X_i,X_j\rangle= 0,\quad \langle x^i,x^j\rangle=0, \quad
\langle x^i,X_j\rangle= \delta^i_j,\quad \forall i,j.\label{age}
\ee
This pairing is essential when Lie bialgebras are connected with
Poisson-Lie groups. Then, the $x^i$ generators are interpreted as
local coordinates on the Lie group with Lie algebra $\G$, and the tensor
$f^{ij}_k$ gives a Poisson bracket among this coordinates. Such a
bracket is just the linear part of the full Poisson-Lie structure
on the Lie group associated to Lie bialgebra $(\G,\eta)$.

\smallskip
\noindent $\bullet$ {\it Example 1.} Let us now explicitly consider 
the Euclidean Lie   algebra $e(2)$ with Lie brackets
\be
[J_{12},P_1]=P_2,\qquad [J_{12},P_2]=-P_1,\qquad [P_1,P_2]=0.
\label{euc}
\ee
The cocommutator
\be
\eta(J_{12})= z\,J_{12}\wedge P_2,\qquad
\eta(P_1)=z\,P_1\wedge P_2,\qquad
\eta(P_2)= 0,
\label{exa}
\ee
endows $e(2)$ with a one-parameter family of
 Lie bialgebra structures. Let us write
explicitly  the Double Lie algebra linked to it by writing the generators
of $e(2)^\ast$ as $\{j_{12},p_1,p_2\}$:
\bea
&& [J_{12},P_1]=P_2,\qquad\qquad [J_{12},P_2]=-P_1,
\qquad\qquad [P_1,P_2]=0,\cr
&& [j_{12},p_1]=0,\qquad\quad\qquad [j_{12},p_2]=z\,j_{12},
\qquad\qquad [p_1,p_2]=z\,p_1,\cr
&& [j_{12},J_{12}]=-z\,P_2,\qquad\quad [j_{12},P_1]=0,
\qquad\quad\qquad [j_{12},P_2]=0,\cr
&& [p_1,J_{12}]=-p_2,\qquad\qquad [p_{1},P_1]=-z\,P_2,
\qquad\qquad [p_{1},P_2]=j_{12},\cr
&& [p_2,J_{12}]=p_1 + J\,p_{12},\quad 
[p_{2},P_1]=-j_{12}+ z\,P_1,\qquad [p_{2},P_2]=0.
\label{doue2}
\eea
Note that we have assumed a commutator notation for
the Lie bracket in the Double. We have also explicitly preserved the
parameter $z$ within the original cocommutator. Such a a parameter will
be identified with the deformation parameter of the
Euclidean quantum algebra that has (\ref{exa}) as the first
order deformation of the quantum coproduct. In this context, the
explicit consideration of $z$ will appear as relevant when contractions
at the quantum algebra level are discussed (see \cite{BGHOS}).

As we have commented before, the six dimensional algebra (\ref{doue2})
admits, by construction, a coboundary Lie bialgebra structure $\delta$
that can be deduced from the r-matrix (\ref{rmat}) (in this case,
$r=J_{12}\otimes j_{12} +
P_1\otimes p_1 + P_2\otimes p_2$):
\bea
&& \delta(J_{12})= z\,J_{12}\wedge P_2,\qquad
\delta(P_1)=z\,P_1\wedge P_2,\qquad
\delta(P_2)= 0,\cr
&& \delta(j_{12})= 0,\qquad\qquad
\delta(p_1)=\,p_2\wedge j_{12},\qquad
\delta(p_2)=j_{12}\wedge p_1 .
\label{doco}
\eea

\sect{Lie algebra and Lie bialgebra contractions}

We present in this section a brief review about the basic ideas related
with the contraction of Lie algebras and the extension of this procedure
to the case of Lie bialgebra structures.

\subsect{Lie algebra contractions}

Let us start by remembering the concept of contraction of Lie algebras
according with the Saletan approach \cite{Sal}.

Let $(A,m)$ and $(A',m')$ be two algebras with the same underlying vector
space $V$ and products $m$ and $m'$, i.e., $m:V\otimes V\to V$. 
Let us assume that there exists a continuous uniparametric family
$\phi_\varepsilon$ of linear mappings
\be
\phi_\varepsilon:V\rightarrow V,\qquad \varepsilon\in (0,1],\qquad
\phi_\varepsilon\big |_{\varepsilon=1}=id,\label{aa}
\ee
such that $\phi_\varepsilon$ is invertible when $\varepsilon\neq 0$ and
singular when $\varepsilon = 0$. The algebra $(A',m')$ is said to
be a contraction of $(A,m)$ if $m'$ can be defined as
\be
m'=\lim_{\varepsilon\to 0}\, m_\varepsilon=\lim_{\varepsilon\to
0}\, \phi_\varepsilon^{-1} \circ m\circ (\phi_\varepsilon\otimes
\phi_\varepsilon).\label{ab}
\ee

When $A$ and $A'$ are Lie algebras with Lie brackets $m$ and $m'$,
the expression (\ref{ab}) can be written as
\be
[X, Y]':=\lim_{\varepsilon\to 0} {\phi_\varepsilon^{-1}
[\phi_\varepsilon(X),\phi_\varepsilon(Y)]}.\label{ac}
\ee

In the case that $V=H\oplus W$, such that $H$ is the underlying
vector space of a subalgebra of $A$ and $W$ is the suplement of $H$, and 
the mapping $\phi_\varepsilon$ is defined by 
\be
\phi_\varepsilon|_H=id, \qquad \phi_\varepsilon|_W=\varepsilon\ id. 
\label{acc}
\ee
we shall say that we have performed an In\"on\"u--Wigner \cite{IW}
contraction along the subalgebra $H$.  

An interesting generalization of the IW contraction can be defined as
follows. Let $\g$ be a Lie algebra whose associated 
vector space $V$ is written as a direct sum of vector subspaces 
\be
V=\bigoplus_i\,V_i,\qquad i=0,1,\dots,N\geq 1. \label{ad}
\ee
The mapping $\phi_\varepsilon$ will be called 
generalized IW contraction \cite{WW} if
\be
\phi_\varepsilon |_{V_i}=\varepsilon^{n_i} \,  Id|_{V_i},\qquad
0\leq n_0 < n_1<n_2<\dots<n_N, \qquad n_i\in\R.\label{ae}
\ee
It can be shown that a given Lie algebra admits a generalized IW
contraction  if and only if
\be
[V_i,V_j]\subset \bigoplus_k V_k,
\label{af}
\ee
where by (\ref{af}) we understand that a given subspace $V_k$ can
originate a contribution to the right hand side of the bracket if 
$n_k\leq n_i + n_j$.

Under a generalized IW contraction the structure constants of the
contracted Lie algebra can  be written in terms of the original ones as
\be
c^{\prime\ k}_{\ ij}= \lim_{\varepsilon\to 0}
\varepsilon^{n_i+n_j-n_k}c_{ij}^k,\label{afa}
\ee
since $\phi_\varepsilon$ acts on any generator $X$ in  $V_i$ as
$\phi_\varepsilon(X)= \varepsilon^{n_i}\,X$.

A further generalization of the above kind of contractions can be
introduced by considering negative values of the $n_i$ exponents
\cite{DM}. Finally, we recall that
the so-called graded contractions, that were introduced in the early
nineties
 by Montigny, Moody
and Patera \cite{dMP}, extended these concepts to  include arbitrary
modifications of the structure constants  of a given algebra compatible
with a fixed grading.

\subsect{Lie bialgebra contractions}

Let $(\g,\eta)$ be a Lie bialgebra and let $\g'$ be a Lie algebra 
obtained from $\g$ by means of a (generalized) IW contraction
 $\phi_\varepsilon$. If $n$ is a positive real number such that the limit
\be
\eta':=\lim_{\varepsilon \to 0}{\varepsilon^n
(\phi_\varepsilon^{-1}\otimes\phi_\varepsilon^{-1})\circ\eta\circ
\phi_\varepsilon} \label{ah}
\ee
exists, then is possible to prove that $(\g',\eta')$ is a Lie bialgebra
\cite{BGHOS}. In fact, there exists a unique minimal value $f_0$ of $n$
such that, if $n\geq f_0$ the limit (\ref{ah}) exists, and if $n>f_0$ that 
 limit is zero. 

We will say that $(\g',\eta')$ is a contracted Lie bialgebra of 
$(\g,\eta)$. The pair $(\phi_\varepsilon,n)$ will be called a  Lie
bialgebra contraction (or bicontraction). The minimal value $f_0$ that
ensures the existence of (\ref{ah}) will be called the fundamental
contraction constant of $(\cal G,\eta)$ associated to the contraction 
$\phi_\varepsilon$, and $(\phi_\varepsilon,f_0)$ will be called
a fundamental bicontraction.

These results can be extended to the case of coboundary Lie
bialgebras \cite{BGHOS}, i.e., the contraction of the
classical $r$--matrix can be considered. If $(\g,\eta(\rho))$ is a
coboundary Lie bialgebra with $\rho$ its classical $r$--matrix and
$\g'$ is a Lie algebra obtained by means of a (generalized) IW
contraction  $\phi_\varepsilon$ from $\g$, we can study the limit
\be
\rho':=\lim_{\varepsilon \to 0}{\varepsilon ^n
(\phi_\varepsilon ^{-1}\otimes \phi_\varepsilon ^{-1}) (\rho)}.\label{ai}
\ee
One can proves that if there is a positive real number $n$ such that the
limit (\ref{ai}) exists, then $(\g',\eta'(\rho'))$ is a coboundary Lie
bialgebra. Also in this case, there is a unique minimal value $c_0$ of $n$
such that, if $n\geq c_0$ (\ref{ai}) exists and, if $n>c_0$ such a limit
is zero.

The minimal value $c_0$ that ensures the existence of the limit (\ref{ai})
will be called coboundary contraction constant of the Lie bialgebra
$(\g,\eta(\rho))$ relative to the contraction mapping
$\phi_\varepsilon$. The pair $(\phi_\varepsilon,c_0)$ is called a
coboundary bicontraction of $(\g,\eta(\rho))$ associated to
$\phi_\varepsilon$. Moreover, the relation $f_0\leq c_0$ is always
fulfilled for a given 
$(\g,\eta(\rho))$ and $\phi_\varepsilon$. 

Complete proofs of all the statements included in this section can be
found in
\cite{BT} as well as various examples (see also \cite{BGHOS,BO}).

\sect{Contracting the classical double}

Since $\cal D$ is a Lie algebra, the theory developed in
subsection 2.1 can be fully applied onto $\cal D$ in order 
to obtain new Lie algebras of the same dimension. However, it would be
interesting to know how it is possible to characterize the contractions
on $\cal D$ that give rise to a double Lie algebra.

\subsect{Double-preserving contractions of $\cal D$}

We shall say that a contraction of $\cal D$ is a ``double-preserving"
contraction, if the contracted algebra ${\cal D}'$ is a classical double,
i.e., if it preserves the internal structure given by (\ref{agd}).

Let us consider the most arbitrary generalized Doebner-Melsheimer
 contraction constructed by considering 
that each  subspace  is spanned by only one generator of $\cal D$. 
This means that the  contraction mapping  is defined as 
\be
\phi_\varepsilon (X_i)=\varepsilon ^{m_i}X_i, \quad
\phi_\varepsilon (x^i)=\varepsilon ^{n_i}(x^i), 
\label{ba} 
\ee 
where $m_i,n_i\in Z$. If we use now the definition (\ref{ac}) to contract
the Double  commutation rules (\ref{agd}), we shall obtain the following
contracted structure tensors $c'$ and $f'$: 
\be\ba  {cc}
c^{\prime\ k}_{\ ij}=\lim_{\varepsilon\to 0}
 {\varepsilon^{m_i+m_j-m_k}c_{ij}^k}, \qquad
f^{\prime\ ij}_{k}=\lim_{\varepsilon\to 0}
  {\varepsilon^{n_i+n_j-n_k}f^{ij}_{k}}\\
\ \ c^{\prime\ k}_{\ ij}=\lim_{\varepsilon\to 0} 
 {\varepsilon^{m_i-n_j+n_k}c_{ij}^k}, \qquad  \
f^{\prime \ ij}_{k}=\lim_{\varepsilon\to 0}
  {\varepsilon^{n_i-m_j+m_k}f^{ij}_{k}}.
\label{bb}
\ea\ee
The expressions in the first row come from the contraction of commutators
$[X_i,X_j]$ and $[x^i,x^j]$, respectively; and those of the second row
from
$[x^i,X_j]$. If we impose the $c'$ and $f'$ components to be the same
whenever they appear in the contracted Double, we are lead to  the
conditions
\be 
m_i-m_j=-(n_i-n_j),
\label{bc}
\ee
for all couples $(i,j)$ such that there exists
either a non vanishing $c_{ij}^k$ or $f^{ij}_k$ component in the original
Double (when a given component of any of the original tensors
vanishes, there is no consistency constrain in the contracted limit
of such a component, that will be zero in any  case -see (\ref{bb})-). 

Afterwards, in order to avoid divergences in the limit $\varepsilon \to 0$
of (\ref{bb}) we shall have to impose that, for all triads $(i,j,k)$
labelling a  non-vanishing component either of  $c$ or of the tensor $f$
\be 
m_i+m_j-m_k\geq 0,\qquad n_i+n_j-n_k\geq 0.
\label{bba}
\ee  
This is tantamount to say that we need the generalized IW contraction
(\ref{ba}) to be a good contraction for both $\g$ and $\g^\ast$ Lie
algebras. Conditions (\ref{bc}) and (\ref{bba}) will therefore suffice to
guarantee that ${\cal D}'$ is a classical double.

\smallskip

\noindent $\bullet$ {\it Example 2.} A quite simple Double-preserving
contraction is obtained if we define $m_i=n_i=N\in Z^+,\ \forall i$. This
contraction always fulfills (\ref{bc}) and (\ref{bba}) and
 originates an abelian
contracted algebra since all the structure constants vanish after
the limit $\varepsilon\to 0$.

\smallskip

\noindent $\bullet$ {\it Example 3.} Let us now consider the 
two dimensional Euclidean double introduced in Example 1 and let us try
to implement on it the non-relativistic limit. It is well
known 
that, at the Lie algebra level,
 such a limit is equivalent to the following IW contraction
\be
\phi_\varepsilon (J_{12})=\varepsilon\,J_{12},\qquad
\phi_\varepsilon (P_1)=P_1,\qquad
\phi_\varepsilon (P_2)=\varepsilon\,P_2,
\label{nonrel}
\ee
where $\varepsilon=1/c$, with $c$ the speed of light. We have to find now
how the contraction mapping $\phi$ has to act on the remaining 
$x^i$ generators of $\cal D$ 
in order to obtain a correct non-relativistic limit on the double.

If we consider
that $\langle J_{12}\rangle=V_1$,$\langle P_{1}\rangle=V_2$ and 
$\langle P_{2}\rangle=V_3$, we shall have that , from (\ref{ba}),
$m_1=m_3=1$ and $m_2=0$.
This notation immediately implies the definition of the remaining
subspaces as
$\langle j_{12}\rangle=v_1$,$\langle p_{1}\rangle=v_2$ and 
$\langle p_{2}\rangle=v_3$. Now we can say that the structure tensors
on the Euclidean double only have as the only non-vanishing components
$c_{12}^3,c_{13}^2,f^{13}_1,f^{23}_2$ and their skew-symmetric
counterparts. This fact leads us to impose onto
the $n_i$ exponents the conditions (\ref{bc}):
\bea
&& m_1-m_2=n_2-n_1=1,\cr    
&& m_1-m_3=n_3-n_1=0,\cr    
&& m_2-m_3=n_3-n_2=-1.    
\label{ecuac}
\eea
The obvious set of solutions for these equations is
\be
n_1=\alpha,\qquad n_2=1+\alpha,\qquad n_3=\alpha, \quad \alpha\in Z.
\label{sol}
\ee
Now we have to go to the remaining conditions (\ref{bba}), that are
translated into the following inequalities, each of them coming from
a triad $(i,j,k)$ labelling  a not-vanishing tensor component:
\bea
&& m_1+m_2-m_3\geq 0,\qquad m_1+m_3-m_2\geq 0,\cr
&& n_1+n_3-n_1\geq 0,\qquad\quad  n_2+n_3-n_2\geq 0.
\label{ineq}
\eea 
The two first ones are obviously fulfilled by the initial values of $m_i$
(we started from a right contraction of the Euclidean Lie algebra).
 The second ones are indeed
equivalent to write $n_3=\alpha\geq 0$. Therefore only generalized
IW contractions are allowed in this case. Moreover, it can be easily
checked that if $\alpha > 0$, all the components of $f'$ vanish, thus
obtaining  a Double algebra of a (1+1) Galilei bialgebra with trivial
cocommutator. On the other hand, if we compute the full contracted
algebra when $\alpha=0$, i.e., with contraction mapping given by
(\ref{nonrel}) and
\be
\phi_\varepsilon (j_{12})=j_{12},\qquad
\phi_\varepsilon (p_1)=\varepsilon\, p_1,\qquad
\phi_\varepsilon (p_2)=p_2,
\label{nonreldu}
\ee
we obtain a expresion for a classical Galilei double that differs with
respect to (\ref{doue2}) in the following commutation rules: 
\be
[J_{12},P_2]=0,\qquad [p_1,J_{12}]=0,\qquad [p_1,P_2]=0.
\label{dogal}
\ee
In fact, we have just obtained the double of 
a Galilei Lie bialgebra that preserves under contraction
the same cocommutator (\ref{exa}). As we have proven by the  previuos
analysis, this is the double corresponding to the
 only non-trivial Galilei Lie
bialgebra that can be reached from our Euclidean one by using the
non-relativistic limit.

\subsect{The duality problem}

An important property of the double is the pairing
$\langle x^i,X_j\rangle= \delta^i_j$. In what follows we shall analyse
the implications of a duality preservation during the contraction 
process.

Let $\phi_\varepsilon$ be the contraction mapping defined by (\ref{ba}).
Then  we shall say that
the pairing $\langle x^i,X_j\rangle$ will be preserved under contraction
if the ``contracted generators"
$\phi_\varepsilon(x^i),\phi_\varepsilon(X_j)$  fulfill
\be
\langle\phi_\varepsilon(x^i),\phi_\varepsilon(X_j)\rangle=
\varepsilon ^{n_i+m_j}\langle x^i,X_j\rangle=\delta^i_j, \label {bd}
\ee
therefore, we get the condition
\be
n_i+m_i=0. \label{be}
\ee
And this means that we are forced to use negative  powers in order to
define contractions preserving the  pairing. We thus deal now with
Doebner-Melsheimer contractions \cite{DM}. 

Let us now suppose that we start from a known contraction of the original
Lie algebra $\cal g$ that is a generalized IW one. This implies that
all $m_i>0$ and, therefore, all $n_i<0$. This fact seems to indicate that
we can find serious problems in the convergency of the $\varepsilon\to  0$
limit (or, equivalently, in solving conditions (\ref{bba})). 

These difficulties
can be confirmed by making use of the Euclidean 
example we have worked out previously and considering
the same contraction mapping (\ref{nonrel}) on the generators of the
$e(2)$  algebra.
If we want to preserve the pairing (\ref{age}) under contraction, we are
forced to have $n_1=-1$, $n_2=0$ and $n_3=-1$. Thus, $\alpha$ has to be
$-1$, a value which is forbidden by convergency conditions (\ref{ineq}).
In conclusion, for this classical double, and provided the contraction
mapping for the Euclidean Lie algebra generators is given by
(\ref{nonrel}),
 there does not exist any
contraction (in the sense of the previuos section)
that preservs the pairing (in the sense of (\ref{bd})).

A possibility in order to 
avoid such divergencies in the $f'$ components 
is to use a ``renormalization" factor 
$\varepsilon ^{t_0},\ t_0\geq 0 \in \R$, such that $\forall N\geq t_0$,
and for all the triads $(i,j,k)$ coming from non-vanishing
components of $f$, the new definition of the 
contracted components be
\be
{f'}^{ij}_k:=\lim_{\varepsilon\to 0}
{\varepsilon^{n_i+n_j-n_k+N}f^{ij}_{k}}.
\label{newf}
\ee
By construction, this contracted tensor does not diverge (of 
course, the exponent $t_0$ depends on both the
mapping $\phi_\varepsilon$ and the double considered). 

However, this way for contracting the double  and preserving
duality has to be shown to be consistent
 (i.e., we should prove now that definition (\ref{newf}) and
the usual contraction of the tensor $c$ give rise to a Lie algebra
with classical double structure). Fortunately, further
computations are not necessary
 because (\ref{newf}) is nothing but the rephrasing in the
classical double language of  the Lie bialgebra  contraction theory 
of Section 2. It can be easily  checked  that
$t_0$ is just the fundamental contraction constant $f_0$, and
that the right ``renormalized" contraction gives always rise to 
the classical double
corresponding to the contracted Lie bialgebra obtained by making use of
the theory sketched in Section 3.2.

In particular, if this duality-preserving approach is applied onto
the Euclidean double (\ref{doue2}), we shall obtain that $t_0=1$. 
On the other hand, the Euclidean Lie bialgebra  (\ref{euc}-\ref{exa})
can be contracted by means of a fundamental
bicontraction characterized, as expected, by $f_0=1$. Moreover, the
Galilean double so obtained coincides with the one derived
in the previuos section by using a generalized IW contraction. 

\subsect{The classical limit as a contraction}

Given a quantum algebra, the limit $z\to 0$ of the deformation parameter
can be interpreted as the classical limit: under it, the quantum
coproduct in the deformed algebra 
becomes cocommutative and, consequently, the
algebra of functions on the group becomes commutative. At the classical
double level, this procedure is equivalent to make zero all the 
components of the tensor $f$ \cite{BGd}. 
These assertions can be easily illustrated by using
the Example 2.

It is remarkable that such a limit  process can be described in a
completely equivalent way as a pure IW contraction on the classical 
double along the Lie (sub)algebra $\g$.
 Explicitly, let us consider an arbitrary classical double and
the contraction mapping  on it  defined by
\be
\phi_\varepsilon (X_i)=X_i, \quad
\phi_\varepsilon (x^i)=\varepsilon\,x^i. 
\label{bacl} 
\ee 
Obviuosly, since $m_i=0$ and $n_i=1$ for all $i$, 
this is a double-preserving  contraction fulfilling conditions (\ref{bc})
and (\ref{bba}). It is also inmediate to check that the contracted double
has commutation rules
\be
[X_i,X_j]= c^k_{ij}X_k, \quad  
[x^i,x^j]=0, \quad
[x^i,X_j]= c^i_{jk}x^k.
\label{agdc}
\ee
Note that another double-preserving contraction can be obtained by
defining
\be
\phi_\varepsilon (X_i)=\varepsilon\,X_i, \quad
\phi_\varepsilon (x^i)=x^i. 
\label{baclb} 
\ee 
This process anhiquilates the $c$ tensor, and can be interpreted
as the classical limit on the dual, that gives rise to
a classical double corresponding to a Lie bialgebra structure  on
a commutative  algebra.
 
\section*{Acknowledgments}
This work has been partially supported by projects
DGICYT/PB94-1115 and DGES/PB95--0719
from the Mi\-nis\-terio de Educaci\'on y Ciencia de Espa\~na. 



\begin{thebibliography}{99}

\bibitem{Dr} V.G.Drinfel'd, Sov. Math. Dokl. {\bf 27} (1983), 68.

\bibitem{CP}  V. Chari, A. Pressley, {\em A Guide to Quantum Groups}.
Cambridge Univ. Press, Cambridge 1994.

\bibitem{FRT} N. Yu. Reshetikhin, L.A. Takhtadzhyan and
L.D. Fadeev, Leningrad Math. J. {\bf 1} (1990), 193.

\bibitem{BCGS} E. Celeghini, R. Giachetti, E. Sorace
and  M. Tarlini, {\em Contractions of quantum groups}\/,
Lecture Notes in Mathematics n. 1510.
Springer-Verlag, Berlin (1992).

\bibitem{CGe2} E. Celeghini, R. Giachetti, E. Sorace
and  M. Tarlini, J. Math. Phys. {\bf 32} (1991), 1159.      

\bibitem{Lukap} J. Lukierski, H. Ruegg and A. Nowicky,
Phys. Lett. {\bf B293} (1992), 344.

\bibitem{BHOS4d} A. Ballesteros, F.J. Herranz, 
M.A. del Olmo and  M. Santander, J. Math. Phys. {\bf 35} (1994), 4928.      

\bibitem{BGHOS} A. Ballesteros, N.A. Gromov, F.J. Herranz, 
M.A. del Olmo and  M. Santander, J. Math. Phys. {\bf 36} (1995), 5916.

\bibitem{BT} A. Ballesteros, {\em Contractions of Lie bialgebras
and quantum deformations of kinematical symmetries}\/,
Ph. D. Thesis (in Spanish), Universidad de Valladolid (1995).

\bibitem{BO} A. Ballesteros and M.A. del Olmo, {\em Contractions of 
Poisson-Lie Groups, Lie
Bialgebras and Quantum Deformations}\/, in ``Quantum Groups and Quantum
Spaces", Warsow (1996).

\bibitem{Sal}  E. Saletan, J. Math. Phys. {\bf 2} (1961), 1.

\bibitem{IW}  E. In\"on\"u and E.P. Wigner, Proc. Natl. Acad. Sci. U. S.
{\bf 39} (1953), 510.

\bibitem{WW}  E. Weimar-Woods, J. Math. Phys {\bf 32} (1991), 2028.

\bibitem{DM} H.D. Doebner, O.Melsheimer, Il Nuovo Cim. {\bf A49} (1967),
306.

\bibitem{dMP} M. de Montigny, J. Patera, J. Phys. A: Math. Gen. {\bf 24}
(1991), 525; R.J. Moody ,J. Patera, J. Phys. A: Math. Gen. {\bf 24}
(1991), 2227.

\bibitem{BGd} F. Bonechi, R. Giachetti, E. Sorace
and  M. Tarlini, {\em Quantum double and  differential 
calculi}, Lett. Math. Phys. {\bf 37} (1996), 405.      

\end{thebibliography}
\end{document}